\documentclass[a4paper,14pt]{article}[2000/05/19]
\usepackage[english]{babel}
\usepackage{amssymb,amsthm,amsmath}
\usepackage{inputenc}
\newtheorem{theorem}{Theorem}[section]
\newtheorem{lemma}[theorem]{Lemma}

\newtheorem{corollary}[theorem]{Corollary}
\theoremstyle{definition}
\newtheorem{definition}[theorem]{Definition}

\theoremstyle{remark}
\newtheorem{remark}[theorem]{Remark}
\numberwithin{equation}{section}
\begin{document}
\large
\setcounter{page}{1}

\begin{center}
United Nations Educational, Scientific and Cultural Organization \\
and\\
International Atomic Energy Agency \\
THE ABDUS SALAM INTERNATIONAL CENTRE FOR THEORETICAL PHYSICS
\end{center}

\vspace{0.5 cm}

\begin{center}
\textbf{FIXED POINTS OF REFLECTIONS OF COMPACT CONVEX SETS AND A CHARACTERIZATION
OF STATE SPACES OF JORDAN BANACH ALGEBRAS}
\end{center}

\vspace{0.5 cm}

\begin{center}
Sh. A. Ayupov \footnote{Senior Associate of ICTP. Corresponding  author. sh\_ayupov@mail.ru}\\
\textit{Institute of Mathematics and Information  Technologies, Uzbekistan Academy of Sciences
Dormon yoli str., 29, 100125,  Tashkent,   Uzbekistan} \\
\textit{and} \\
\textit{The Abdus Salam International Centre for Theoretical Physics, Trieste, Italy}
\end{center}

\begin{center}
and
\end{center}

\begin{center}
N.~J.~Yadgorov\\
\textit{Institute of Mathematics and Information  Technologies, Uzbekistan Academy of Sciences
Dormon yoli str., 29, 100125,  Tashkent,   Uzbekistan}
\end{center}

\vspace{0.5 cm}

\begin{abstract}
In the present article we prove a fixed point theorem for reflections of compact convex sets and give a new characterization of state space of \textit{JB}-algebras among compact convex sets. Namely they are exactly those compact convex sets which are strongly spectral and symmetric

\end{abstract}

\vspace{1 cm}

\begin{center}
MIRAMARE --- TRIESTE
\end{center}

\newpage

\section{Introduction}

\hspace{0.4cm} The present paper is devoted to the operational or "convex" approach to the axiomatics of quantum theory, the basic concept of which is the convex set of states of a physical system.

In the well-known algebraic approach to quantum field theory the space of observables from a Jordan Banach algebra (\textit{JB}-algebra), i.e.
a Jordan algebra $A$ over reals with an identity element $e$ equipped with a complete norm such that for $a,b, \in A:$
$$\lVert a^2\rVert=\lVert a\rVert^2, \quad \lVert a^2\rVert\leq \lVert a^2+b^2\rVert.$$

Recall that if $A$ is a \textit{JB}-algebra, then the set $A^+$ of all the squares in $A$ is a proper convex cone organizing $A$ to a (norm) complete order-unit space, whose distinguished order-unit is the multiplicative identity $e,$ whose norm is the given one, and such that for $a\in A$
\begin{flushright}
$-e\leq a\leq e$ implies $0\leq a^2\leq e.$ \hspace{4cm} $(\ast)$
\end{flushright}

Conversely, if $A$ is a complete order-unit space, equipped with a Jordan product, for which the distinguished order-unit acts as the identity element and such that $(\ast)$ is satisfied, then $A$ is a \textit{JB}-algebra in the order-unit norm.

As a rule, we shall consider \textit{JB}-algebras which are Banach dual spaces, i.e. $A=V^{\ast}$ for some Banach space $V,$  and we will refer to them as \textit{JBW-algebras.} Probability measures in this algebraic approach correspond to states on $A$ (normal states as usual). Recall that a \textit{state} is a positive linear functional $\rho$ on $A$ such that $\rho(e)=1.$ A functional $f$ is said to be normal if $f(a)=\lim f(a_{\alpha}),$ whenever $\{a_\alpha\}$ is an increasing net in $A$ with $\sup a_\alpha=a.$  It is known that if $A$ is a \textit{JBW}-algebra, then its predual $V$ is unique and it can be identified with the space of all normal linear functionals in $A^\ast.$

Let us consider some examples of \textit{JB}- and \textit{JBW}-algebras.

\begin{enumerate}
    \item  The self-adjoint part $W_{sa}$ of a $C^\ast$-algebra (resp. von Neumann algebra) $W$ with the symmetrized product $a\circ b=(ab+ba)$ is a \textit{JB}-algebra (resp. \textit{JBW}-algebra).
    \item The algebra $L_\mathbb{R}^\infty(\Omega, \mu)$ of all bounded random variables on a classical probability space $(\Omega, \mu)$ is an associative \textit{JBW}-algebra.
    \item The exceptional Jordan algebra $M_3^8$ of all symmetric $3\times 3$ matrices over the Cayley numbers is a finite dimensional \textit{JB}-algebra (and therefore a \textit{JBW}-algebra).
    \item \textit{Spin factors} can  be defined as follows. Let $H$ be a real Hilbert space. Consider the vector space $R\times H$ of pairs $(\alpha, h), \alpha\in R, h\in H,$ with the product
    $$(\alpha, h)\circ(\beta, g)=(\alpha\beta+\left\langle h, g\right\rangle, \alpha g+\beta h),$$
where $\left\langle h, g\right\rangle$ is the inner product of the vectors $h,g\in H.$ Then $A=R\times H$ becomes a Jordan algebra, which is a \textit{JBW}-algebra with respect to the norm $$\lVert(\alpha, h)\rVert=\lvert\alpha\rvert+\lVert h\rVert.$$

Such \textit{JBW}-algebras are called spin factors, they are exactly \textit{JBW}-factors of the type $I_2$ (for detail, see \cite{16}).
\end{enumerate}
It is well known that the state space of a \textit{JB}-algebra (in particular $C^\ast$-algebra) is a compact convex set (a simplex in the classical associative case). The converse problem is essentially more interesting and difficult: to characterize these state spaces among general compact sets in a locally convex space. That is to find geometric conditions for a convex set $K$ to be affinely isomorphic and homeomorphic to the state spaces of a $C^\ast$-algebra or a \textit{JB}-algebra and, more general, to the normal state spaces of a von Neumann algebra or a \textit{JBW}-algebra.

This problem is interesting in its own right and is very important for applications in the operational (convex)  approach to the axioms of the quantum theory. Various geometric and physical conditions of this kind have been suggested in \cite{2}, \cite{4}, \cite{6}-\cite{8}, \cite{10}-\cite{14}, \cite{17}.

In the present paper we give the most simple geometric conditions for a compact convex set to be the state space of a \textit{JB}-algebra. Namely they are exactly those compact convex sets which are strongly spectral and symmetric. Earlier a similar result have been obtained for the finite dimensional case in \cite{10}, \cite{11}; for the modular \textit{JBW}-algebras in \cite{12} and the semi-modular case in \cite{13}.

\section{Affine function spaces on convex sets}

\hspace{0.4cm}Let $K$ be a convex set in a locally convex space $V$, and let $A=A^b(K)$ denote the space of all bounded affine functions on $K$ with pointwise ordering. Then $(A, e)$ is an order-unit space, where $e=1$ is the distinguished order-unit. Without loss of generality, one can assume that $K$ is regularly imbedded into $V,$ i.e. $(V, K)$ is a base norm space, such that $(A, e)$ an $(V, K)$ are in separating order and norm duality, and $A=V^\ast$ (here and below we refer to \cite{3} and \cite{4} for details).

\begin{definition}
A positive norm one projection $R:A\rightarrow A$ is said to be $P$-projection if there exists a unique positive norm one projection $R':A\rightarrow A$ such that
\begin{align*}
    \text{im}^+R=\ker^+R', & \qquad \text{im}^+R^\ast=\ker^+R'^{\ast},\\
    \ker^+R=\text{im}^+R', & \qquad \ker^+R^\ast=\text{im}^+R'^\ast,
\end{align*}
where $R^\ast$ is the dual projection for $R,$ i.e. $R^\ast:V \rightarrow V$ and $Ra(\rho)=a(R^\ast\rho)$ for $a\in A, \rho\in V.$
\end{definition}
Denote $\mathfrak{B}$ the set of all $P$-projections on $A,$ and define an ordering in it by $R\leq Q$ when $\text{im} R \subseteq \text{im} Q,$ i.e. $RQ=QR=R.$ It is clear that $O\leq R\leq I,$ where $O$ is the zero projection and $I$ is the identity map. The projection $R'$ is also a $P$-projection, which is called the \textit{quasicomplement} for $R.$

To every $P$-projection $R$ is associated a \textit{projective unit} $u=Re\in A.$ Denote by $\mathfrak{U}$ the set of all projective units in $A$ and consider the \textit{natural ordering on $\mathfrak{U},$ induced from $A,$ and the orthocomplementation} : $Re \rightarrow (e-Re).$

Recall that a convex subset $G$ of $K$ is called a face if for $x,y \in K$ and $0<\lambda<1$ the relation $\lambda x+(1-\lambda)y \in G$ implies that $x,y \in G.$ A point $z\in K$ is said to be an extreme point if $\{z\}$ is a face of $K;$ denote by $\partial_eK$ the set of all extreme points of $K.$ A convex set $K$ is said to be \textit{strictly convex} if any proper face of $K$ is an extreme point. A face $G$ of $K$ is \textit{exposed} if $G=\{\rho\in K: a(\rho)=0\}$ for an appropriate $a\in A^+.$ If, in addition, $a$ is a projective unit, i.e. $a=Re\in\mathfrak{U},$ $R$ is a $P$-projection, then $G$ is called a \textit{projective face}. In other words, the projective faces are the faces of the form
$$G=\text{im}^+R^\ast\cap K=F_R, \quad R \in \mathfrak{B}.$$

Denote by $\mathfrak{F}$ the set of all projective faces of $K,$ and consider the natural (set-theoretical) ordering and the orthocomplementation
$$F_R\rightarrow F_R^\#=F_{R'}=\text{im}^+R'^\ast\cap K,$$
where $R'$ is the quasicomplementary $P$-projection for $R.$ The projective face $F_R^\#$ is called the quasicomplement of $F_R.$

\begin{definition} A convex set $K$ is said to be \textit{projective} if its every exposed face is projective.
\end{definition}
Elements $a,b \in A^+$ are said to be orthogonal (denoted $a\bot b$) if there is a projective face $F$ in $K$ such that $a(\rho)=0$ for all $\rho \in F$ and $b(\sigma)=0$ for all $\sigma\in F^\#.$

\begin{definition} A projective convex set $K$ is said to be \textit{spectral} if any element $a\in A =A^b(K)$ can be uniquely decomposed as $a=a_+-a_-$ with $a_+, a_-\in A^+,\quad a_+\bot a_-.$
\end{definition}
\begin{theorem} (Alfsen and Shultz, \cite{3}). If $K$ is a spectral convex set, then $\mathfrak{B,U}$ and $\mathfrak{F}$ are mutually isomorphic complete orthomodular lattices (the quantum logics).
\end{theorem}

Consider some examples of spectral convex sets.
\begin{enumerate}
    \item  The set $K$ of all normal states on an arbitrary \textit{JBW}-algebra $A$ (in particular, on any von Neumann algebra $W$) is a spectral convex set. In this case $A=A^b(K)$ (resp. $W_{sa}=A^b(K)$) and the notions of $P$-projections, projective units and projective faces coincide with the maps  $x \rightarrow U_px$ $(x\to pxp),$ where $p$ is an idempotent (a projection), the notion of idempotents (projections), and the notion of closed faces of the state spaces of the algebra, respectively. It should be noted that if $K=S(A)$ is the space of all states on a \textit{JB}-algebra $A$ (or a $C^\ast$-algebra), then it is $^\ast$-\textit{weakly compact convex set} and also spectral, since it can be identified with the normal state space of the enveloping \textit{JBW}-algebra (resp. von Neumann algebra) $\tilde{A}$ \cite{9}, \cite{16}.
    \item Let $K_1$ be the three dimensional convex set from Fig. 10 in (\cite{3}, Sect. 10, p.94) which combines "simplicial" and "rotund" features in a slightly less trivial way than the cone. This set is thought of as a "compressed ball" with a "triangular equator". It admits a unique tangent plane at each point of the surface except at the vertices of the triangle. This "compressed ball" is also a spectral convex set.
    \item Let $K_2$ be the unit ball of the $L_p$-space with $1<p<\infty.$ Then $K$ is a spectral convex set (\cite{3}, Theorem 10.3).
    \item If $K_3$ is the set of all points $(x; \,\, y_1, \ldots, y_m; \,\, z_1, \ldots, z_n;\,\, t)\in \mathbb{R}^{m+n+2}$ where $m,n=0, 1, 2, \ldots$ (if $m=0$ or $n=0,$ then there are no $y$-terms or $z$-terms) which satisfy the inequality $t^4\leq \left(x^2-\sum_{i=1}^{m}{y_i^2}\right)\left((1-x)^2-\sum_{j=1}^{n}z_j^2\right),$ together with the inequalities $0\leq x \leq 1,$ $\sum_{i=1}^m{y_i^2}\leq x^2,$ $\sum_{j=1}^n{z_j^2}\leq (1-x^2),$ then $K_3$ is a non-decomposable spectral convex set. (see \cite{7}, Theorem 8.87).

\end{enumerate}

\begin{remark}\label{R1}
The sets $K_1$ , $K_2$ $(p \neq 2)$  and $K_3$ $(m+n \neq 0, m,n = 0,1,2,...)$ are the examples of spectral convex sets which are not affinely isomorphic to the normal state of any \textit{JBW}-algebra.
\end{remark}

\begin{definition} A spectral convex set $K$ is said to be \textit{strongly spectral} if any $a\in A(K),$ where $A(K)$ is the space of all continuous affine functions on $K,$ can be uniquely decomposed as $a=a_+-a_-$ with $a_+, a_-\in A(K)^+,$ $a_+\bot a_-.$
\end{definition}

\textbf{\textbf{
\section{A fixed point theorem for reflections of compact convex sets}}}

Let $K$ be a compact convex set in a locally convex Hausdorff space $V,$ and denote by $\Gamma(K)$ the group of all affine homeomorphisms of $K$ onto itself. For $T \in \Gamma(K)$ (respectively for $G\subset\Gamma(K)$) consider the set
$$K_T=\{p\in K: T(p)=p\}$$
respectively
$$K_G=\bigcap_{T\in G} K_T,$$
of all fixed points of the map $T,$ respectively of all common fixed points of the family $G.$

Recall, that given a compact convex set $K,$ an affine homeomorphism $T:K\rightarrow K$ is called \textit{a reflection}, if $T^2=\textit{id}$ - the identical map. The set of all reflections of the set $K$ is denoted by $S(K),$ i.e.
$$S(K)=\{T\in \Gamma(K): T^2=\text{id}\}$$
Given a subset $K$ in the vector space $V,$ the dimension $\text{dim}K$ means the dimension   of its affine span $\textit{affK}$ , i.e.
$$\text{dim}K:=\text{dim }(\textit{affK}).$$
Recall the following well-known results.
\begin{lemma}\label{2.1} (\cite{18}, page 152. lemma) Let $K$ be a compact convex set in a locally convex Hausdorff space $V,$ and let $T:K\rightarrow K$ be an affine map. Then $T$ has a fixed point.
\end{lemma}

\begin{lemma}\label{2.2} (\cite{15}, page 498, lemma 2.2) Let K be a finite dimensional compact convex set in a locally convex Hausdorff space $V.$ Then the group $\Gamma(K)$ of all affine homeomorphisms of $K$ onto $K$ has a common fixed point.
\end{lemma}

Now we shall prove the following auxiliary results for reflections.
\begin{lemma}\label{2.3} Let $K$ be a compact convex set in a locally convex Hausdorff space $V.$ The any two reflections $T$ and $S$ of $K$ has at least one common fixed point.
\end{lemma}

\textbf{Proof.} Since $T$ and $S$ are reflections, $S\circ T$ is an affine homeomorphism of $K$ onto itself. By Lemma \ref{2.1} it has a fixed point say $p \in K,$ i.e. $S\circ T(p)=p.$ Therefore $S\circ S\circ T(p)=S(p)$ and since $S\circ S= \textit{id},$ we obtain that $T(p)=S(p).$  Consider the point $p_0=\frac{1}{2}(p+T(p))=\frac{1}{2}(p+S(p)).$ Then it is clear that $T(p_0)=S(p_0)=p_0.$ The proof is complete $\square$

Now we shall give the main result of this section.
\begin{theorem} Let $K$ be a compact convex set in locally convex Hausdorff space $V.$  Consider the family $G=\{T\in S(K): \text{dim}K_T<\infty\}$ of reflections with finite dimensional sets of fixed points. Then $G$ has a common fixed point.
\end{theorem}

\begin{proof} Consider $S_1, \ldots, S_n\in G.$ By Lemma \ref{2.1} each $K_{S_i}$ is a non empty compact convex set. Put $E =  \textit{aff} (\cup_{i=1}^n K_{S_i})$ -- the affine subspace of $V$ generated by the sets $K_{S_1}, \ldots, K_{S_n}.$  It is clear that $\text{dim}E<\infty.$  Denote $K_0=E\cap K$ and let us show that $S_i(K_0)\subset K_0$ for each $i=1, \ldots, n.$ From $S_i^2= \textit{id}$ it follows that if $p \in K,$ then $\frac{p+S_i(p)}{2}=p^0\in K_{S_i},$ i.e. $S_i(p)=2p^0-p\in K.$

Let $p\in K_0,$ i.e. $p=\lambda_1p_1+\cdots +\lambda_np_n,$ where $p_i \in K_{S_i}, \lambda_i\in \mathbb{R}, \lambda_1+\cdots+\lambda_n=1.$ We have
$$
S_i(p)=\sum_{i=1}^nS_i(\lambda_ip_i)=\sum_{i=1}^n\lambda_iS_i(p_i)
=\sum_{i=1}^n\lambda_i(2p_i^0-p_i)=2\sum_{i=1}^n\lambda_ip_i^0-\sum_{i=1}^n\lambda_ip_i.
$$

Since $\sum_{i=1}^n\lambda_ip_i^0\in E$ and $\sum_{i=1}^n\lambda_ip_i=p\in E,$ it follows that $S_i(p)\in E,$ i.e. $S_i(E)\subset E.$ Therefore $S_i(K_0)\subset K_0.$ From $\text{dim}E<+\infty$ it follows that $K_0$ is compact and by Lemma \ref{2.1} each $S_i$ has a fixed point in $K_0$ for $i=1, \ldots, n.$

Since $S_1, \ldots, S_n\in G \subset\Gamma(K_0)$ - the group of all affine homeomorphisms of $K_0$ onto itself Lemma \ref{2.2} implies that there exists a common fixed point $p \in K_0$ for all $T\in \Gamma(K_0),$ i.e. $T(p)=p$ for all $T\in \Gamma(K_0).$ Therefore $K_{S_1}\cap \ldots \cap K_{S_n}\neq \emptyset$ for any finite family $S_1, \ldots, S_n\in G$ i.e. $\{K_T: T\in G\}$ is a centered family of closed  subsets of $K.$ From compactness of $K$ it follows that $K_G=\cap_{T\in G}K_T\neq \emptyset,$ i.e. $G$ has a common fixed point in $K.$ The proof is complete.
\end{proof}

\begin{remark}\label{2R} Let $K$ be the state space of a \textit{JBW}-factor $M$ of type $I_\infty.$ Then $S(K)=\{S^\ast: S=2U_p+2U_{p'}-I, p \text{\,\,- a projection in} \,\, M\},$ where $U_p:M\rightarrow M$ is defined as $U_p(x) = 2p(px)-px, x\in M,$ and $I = \textit {id}$. In this case any two such symmetries (reflections) has a common fixed point, but $S(K)$ does not have a common fixed point (a tracial state).
\end{remark}

\section{A characterization of state spaces of Jordan Banach algebras}

\hspace{0.4cm}Let $K$ be a compact convex set in a locally convex Hausdorff space $V.$ As in the Section 2 denote by $A=A^b(K)$ (respectively $A(K)$) the space of all bounded (respectively continuous) affine functions on $K$ with the pointwise ordering. Taking as an order-unit the function $e,$ which is identically equal to 1 on $K$ we obtain that $(A^b(K), e)$ and $(A(K), e)$ are order-unit spaces. Moreover without loss of generality we may assume that $K$ is regularly imbedded into $V$ (see Section 2)

\begin{definition} A convex set $K$ is said to be \textit{symmetric}, if $S_ R=2R+2R'-I\geq 0,$  i.e. $S_R$ is a positive linear operator on the order-unit space $A=A^b(K)$ for each $P$-projection $P\in\mathfrak{R}.$
\end{definition}

By Lemma 3.13 \cite{6} the symmetricity of $K$ means that $K$ is symmetric with respect to the convex hull $co(F\cup F^\#)=F\oplus_cF^\#$ for each projective face $F\in\mathfrak{F},$ i.e. there exists a reflection $T=(2R+2R'-I)^\ast$ on $K$ with set of fixed points exactly equal to $F\oplus_cF^\#$
\begin{remark}\label{2.2R} Let $K$ be a spectral and symmetric compact convex set. Then  every projective face $F$ of $K$ is itself a spectral symmetric set.
\end{remark}

Recall that a convex set $K$ has the \textit{Hilbert ball property} if for any pair $\rho, \sigma$ of extreme points of $K$ the face $face(\rho, \sigma)$ generated by these points, is an exposed face affinely isomorphic to a Hilbert ball, i.e. the closed unit ball of some real Hilbert space.

Now let $L$ be a lattice and $a,b\in L.$ We say that $(a, b)$ is a modular pair (denoted $(a, b)M$) if for all $x \in L$ with $a\wedge b \leq x \leq b$ we have the equality $x=(x\vee a)\wedge b.$ A lattice $L$ is called \textit{semi-modular}, if the relation $M$ is symmetric, i.e. $(a, b)M$ implies $(b,a)M.$ A spectral convex set  $K$ is said to be semi-modular, if the lattice of its projective faces $\mathfrak{F}$ is a semi-modular lattice.

Earlier various characterizations of the state space of \textit{JB}-algebras among compact convex sets have been obtained. Recall two of them
\begin{theorem}\label{T3.1}(\cite{4}, Theorem 7.3) A compact convex set $K$ is affinely and topologically isomorphic to the state space of a \textit{JB}-algebra (with $^\ast$-weak topologically) if and only if $K$ is symmetric, strongly spectral and has the Hilbert ball property.
\end{theorem}

In the following theorem we have replaced the "local" condition "Hilbert ball property" by a "global" condition of semi-modularity.
\begin{theorem}\label{T3.2} (\cite{13}, Theorem 4.4) A compact convex set $K$ is affinely and topologically isomorphic to the state space of a
\textit{JB}-algebra (with the $^\ast$-weak topology) if and only if $K$ is strongly spectral, symmetric and semi-modular.
\end{theorem}

In the paper (\cite{13}, page 8) we have conjectured that in fact both conditions "Hilbert ball property" and "semi-modularity" seem to be redundant. The following main result of the present paper gives the affirmative answer to this conjecture

\begin{theorem}\label{T3.3} A compact convex set $K$ is affinely and topologically isomorphic to the state space of a \textit{JB}-algebra (with the $^\ast$-weak topology) if and only if $K$ is strongly spectral and symmetric.
\end{theorem}

In order to prove this theorem we need several preliminary results
\begin{lemma}\label{3.4}(\cite{12}, Theorem 2.8) Let $K$ be a projective convex set. Then each extreme point of $K$ is a projective face.
\end{lemma}

\begin{lemma}\label{3.5}(\cite{12}, Theorem 3.2; \cite{13},  Theorem 3.4) A strictly convex $K$ is affinely homeomorphic to the state space of a spin factor if and if only if $K$ is strongly spectral and symmetric
\end{lemma}
\begin{lemma}\label{3.5} Let $K$ be a projective compact convex set. Then $K$ is strictly convex if and only if $\mbox{dim}(\{\omega\}\oplus_c\{\omega\}^\#)=1$ for each extreme points $\omega$ in $K.$
\end{lemma}

\begin{proof} Denote by $\partial_{e}K$ the set of all extreme point of $K.$ By Krein-Milman theorem $\partial_{e}K\neq\emptyset.$ Suppose that $K$ is strictly convex. Then by Lemma 4.6 $\{\omega\}$ and $\{\omega\}^\#$ are minimal projective faces of $K$ (here and further we shall identify the extreme point $\omega$ with the minimal projective face $\{\omega\}$). Therefore $\{\omega\}\oplus_c\{\omega\}^\#=[\omega, \omega^\#],$ i.e. $\text{dim}(\{\omega\}\oplus_c\{\omega\}^\#)=1.$

Conversely, suppose that $\text{dim}\left(\{\omega\}\oplus_c\{\omega\}^\#\right)=1$ for each $\omega\in\partial_{e}K.$ Let us show that $K$ is strictly convex. Denote by $\partial K$ the affine boundary of the set $K$ (see for details \cite{12}). If $\tau\in\partial K,$ then the smallest projective face $F(\tau)$ containing $\tau$ is a proper face, i.e. $F(\tau)\neq K$ (\cite{12}, Lemma 2.10).

By Krein-Milman theorem $F(\tau)$ contains an extreme point, say $\omega \in \partial_{e}K,$ i.e.
$\{\omega\}\subset F(\tau),$ and hence $F(\tau)^\#\subset\{\omega\}^\#.$ But $\text{dim}\left(\{\omega\}\oplus_c\{\omega\}^\#\right)=1,$ i.e. $\omega^\#\in \partial_{e}K.$ Therefore $F(\tau)^{\#}=\{\omega\}^{\#},$ i.e. $F(\tau)=\{\omega\}.$ This means that $\tau=\omega \in \partial_{e}K,$ i.e. each point of the affine boundary of $K$ is an extreme point, i.e. $K$ is strictly convex. The proof is complete.

\end{proof}

\begin{lemma}\label{3.6} Let $K$ be a strongly spectral and symmetric compact convex set. Then given arbitrary two extreme points $\rho$ and $\sigma$ in $K$ the projective face $F(\rho, \sigma)$ generated by these points is either strictly convex or coincides with the segment $[\rho, \sigma].$
\end{lemma}

\begin{proof} Since $F(\rho, \sigma)$ is a projective face, it is also strongly spectral and symmetric (see, Remark \ref{2.2R}). So without loss of generality we may assume that $K=F(\rho, \sigma).$ Suppose that $K=F(\rho, \sigma)\neq [\rho, \sigma].$ Thus (by Lemma \ref{3.5})in order to prove the lemma it is sufficient to show that $\text{dim}\left(\{\omega\}\oplus_c\{\omega\}^\#\right)=1$ for each extreme point $\omega$ in $K$. Note also that Proposition 8.86 from \cite{7} implies that if $K_1$ and $K_2$ are spectral sets then $K_1\oplus_cK_2$ is also spectral. Therefore if $K$ is a symmetric spectral set then given any extreme point $\rho \in \partial_cK$ there exists $S_\rho\in S(K)$ such that $K_{S_\rho}=\{\rho\}\oplus_c\{\rho\}^\#$ and since $\{\rho\}$ and $\{\rho\}^\#$ are also  symmetric spectral sets, it follows that $K_{S_\rho}$ is also spectral and symmetric. At the same time $K_{S_\rho}\notin \mathfrak{F},$ i.e. it is not a projective face, if $K_{S_\rho}\neq K.$ Moreover one has that $\partial K_{S_\rho}\subset \partial K$ and $\partial_eK_{S_\rho}=\{\rho\}\cup\partial_e\{\rho\}^\#\subset\partial_eK.$

For each couple $\rho, \sigma$ of extreme points in $K$ there exist $S_\rho, S_\sigma\in S(K)$ such that $K_{S_\rho}=\{\rho\}\oplus_c\{\rho\}^\#,$ $K_{S_\sigma}=\{\sigma\}\oplus_c\{\sigma\}^\#$ (by the symmetricity property and Lemma 3.13 \cite{6}). By Lemma \ref{2.3} we have that $K_{S_\rho}\cap K_{S_\sigma} \neq \emptyset$ and from $F(\rho, \sigma)=K$ it follows that either $\text{dim}(\{\rho\}\oplus_c\{\rho\}^\#)=1$ or $\text{dim}(\{\sigma\}\oplus_c\{\sigma\}^\#)=1.$ Indeed, if both dimensions are more than 1 , then the intersection  $K_{S_\rho}\cap K_{S_\sigma} \neq \emptyset$  is 1-dimensinal and hence contains a point ${\tau}$ from the affine boundary of $K$.  This means that $\rho, \sigma\in F({\tau})\neq K$, in contradiction with the assumption $K=F({\rho, \sigma}).$ So let us suppose without loss of generality that $\text{dim}(\{\rho\}\oplus_c\{\rho\}^\#)=1$. Take an arbitrary extreme point $\omega \in K, \omega\neq\rho, \omega\neq\{\rho\}^\#$. In this case $\{\rho\}\oplus_c\{\rho\}^\# \nsubseteq \{\omega\}\oplus_c\{\omega\}^\#.$ Suppose that  $\text{dim}\left(\{\omega\}\oplus_c\{\omega\}^\#\right)\geq2$.  Then as above by Lemma 3.13 \cite{6} for the extreme point $\omega$ in $K$ there exists a reflection $S_\omega \in S(K)$ such that $K_{S_\omega}=\{\omega\}\oplus_c\{\omega\}^\#$ and by Lemma \ref{2.3} we have $K_{S_\rho}\cap K_{S_\omega} \neq \emptyset.$

Consider $K_1=K_{S_\omega}\cap S_\rho K_{S_\omega}.$ It easy to see that $S_\rho K_{S_\omega} = K_{S_\rho S_\omega S_\rho}$  and that $S_\rho S_\omega S_\rho$  is also a reflection of $K$. Therefore by Lemma \ref{2.3} $K_1\neq \emptyset$.  Since $S_\rho^2=\textit{id}$ it follows that $S_\rho K_1=K_1$ and from $\text{dim} K_{S_\omega}\geq 2$ we have that $\text{dim}{S_\rho}K_\omega \geq 2$ and thus $\text{dim} K_1=\text{dim}(K_{S_\omega}\cap S_\rho K_\omega)\geq 1.$

Further $\text{dim}(\{\rho\}\oplus_c\{\rho\}^\#)=1$ means that $\partial_eK_{S_\rho}=\{\rho, \rho^\#\}$ and since $\rho\notin\{\omega\}\oplus_c\{\omega\}^\#$ it follows that $\partial_eK_{S_\rho}\cap\partial_eK_{S_\omega}=\emptyset$ and $\partial_eK_{S_\rho}\cap\partial_eK_{S_\rho S_\omega S_\rho}=\emptyset.$ Let us prove that $\partial_eK_{S_\omega}\cap\partial_eK_{S_\rho S_\omega S_\rho}=\emptyset.$ Suppose that $\zeta \in \partial_eK_{S_\omega}\cap\partial_eK_{S_\rho S_\omega S_\rho}.$ From $S_\rho K_1=K_1$ it follows that $S_\rho(\zeta)\in K_1\subset K_{S_\omega},$ i.e. both $\zeta$ and  $S_\rho(\zeta)$  belong to $\partial_eK_{S_\omega}.$ Consider the projective face $F(\zeta, S_\rho(\zeta))$ generated by $\zeta$ and $S_\rho(\zeta)$, and note that $F(\zeta, S_\rho(\zeta))\neq K.$ Indeed, since $\zeta, S_\rho(\zeta) \in \partial_eK_{S_\omega},$ $\zeta\neq S_\rho(\zeta),$ and $\text{dim}K_{S_\omega}\geq 2,$ we have the following possibilities:

\textit{a}. $\zeta, S_\rho(\zeta)\in\{\omega\}^\#,$ then it is clear that $F(\zeta, S_\omega(\zeta))\subset\{\omega\}^\#\neq K.$

\textit{b}. $\zeta \notin\{\omega\}^\#$ and $S_\rho(\zeta)\in\{\omega\}^\#$, or respectively, $S_\rho(\zeta)\notin \{\omega\}^\#,$ and $\zeta\in\{\omega\}^\#.$ In this case since $K_{S_\omega}=\{\omega\}\oplus_c\{\omega\}^\#$, it follows that either $\zeta=\omega$ or respectively,  $S_\rho(\zeta)=\omega$. Therefore since $K_{S_\omega}=\{\omega\}\oplus_c\{\omega\}^\#$ and $\text{dim}K_{S_\omega}\geq 2$ it follows that the segment $[\zeta, S_\rho(\zeta)]$ is a subset in $\partial K_{S_\omega}\subset\partial K,$ thus $F(\zeta, S_\rho(\zeta))\subset\partial K$ i.e. $F(\zeta, S_\rho(\zeta))\neq K.$

\textit{c}. If both $\zeta\notin \{\omega\}^\#$ and $S_\rho(\zeta)\notin \{\omega\}^\#,$ then $\zeta=S_\rho(\zeta)=\omega$ that is a contradiction with $\zeta\neq S_\rho(\zeta).$

Therefore in any case $F(\zeta, S_\rho(\zeta))\subset \partial K.$ This implies that $\frac{1}{2}(\zeta + S_\rho(\zeta))\in\partial K$ and it is clear that $\frac{1}{2}(\zeta + S_\rho(\zeta))\in K_{S_\rho}=[\rho, \rho^\#],$  because $S_\rho^2=\textit{id}.$ But since $\rho$ and $\rho^\#$ are extreme points this implies that either $\frac{1}{2}(\zeta+ S_\rho(\zeta))=\rho$ or $\frac{1}{2}(\zeta+ S_\rho(\zeta))=\rho^\#$ and thus either $\zeta=S_\rho(\zeta) =\rho$ or $\zeta=S_\rho(\zeta)=\rho^\#.$ This contradicts to $\partial_eK_{S_\rho}\cap \partial_eK_{S_\omega}=\emptyset.$ Therefore we have proved that $\partial_eK_{S_\omega}\cap \partial_eK_{S_\rho S_\omega S_\rho}=\emptyset,$ in particular $\{\omega\}^\#\cap \{S_\rho(\omega)\}^\#=\emptyset.$

From $\text{dim}K_1 \geq 1$ and $\{\omega\}^\#\cap
\{S_\rho(\omega)\}^\#=\emptyset$ it follows that there exists a
point $\nu \in K_1\cap \partial K$ such that $\omega,
S_\rho(\omega)\in F(\nu)\neq K$ and $F(\omega,
S_\rho(\omega))\subseteq F(\nu).$ By Krein-Milman Theorem there
exists an extreme point $\zeta_0 \in \partial_eF(\nu)^\#\subset
F(\omega, S_\rho(\omega))^\#.$ This implies that $\zeta_0\in
F(\omega, S_\rho(\omega))^\#\subset\{\omega\}^\#$ and $\zeta_0 \in
F(\omega, S_\rho(\omega))^\#\subset\{S_\rho(\omega)\}^\#,$ i.e.
$\zeta_0\in\{\omega\}^\# \cap \{S_\rho(\omega)\}^\#$ is a
contradiction with the above. This contradiction shows that the
assumption $\text{dim}(\{\omega\}\oplus_c\{\omega\}^\#)\geq 2$ is
false, i.e. $\text{dim}(\{\omega\}\oplus_c\{\omega\}^\#)=1$ for
each extreme point $\omega \in \partial_eK.$ The proof is
complete.

\end{proof}

\begin{corollary}\label{c3.7} If $K$ is strongly spectral and symmetric then given arbitrary extreme points $\rho, \sigma$ in $K,$ the face generated by these two points is a projective face, i.e. face $(\rho, \sigma)=F(\rho, \sigma).$
\end{corollary}

\begin{proof}  If $\rho=\sigma$ then $face(\rho, \sigma)=face(\sigma)=\{\sigma\}$ is a projective face by Lemma \ref{3.4}. Thus suppose that $\rho\neq \sigma.$ It is clear that $face(\rho, \sigma)\subset F(\rho, \sigma).$ By Lemma \ref{3.6} $F(\rho, \sigma)$ is strictly convex and $face(\rho, \sigma)$ is a face of $F(\rho, \sigma).$ This is possible only if $face(\rho, \sigma)=F(\rho, \sigma).$ The proof is complete.
\end{proof}

\textit{Proof of Theorem \ref{T3.3}} The necessity is clear and follows from the above Theorem \ref{T3.1}.

Sufficiency, By Lemma \ref{3.6} and Corollary \ref{c3.7}  the face
$face(\rho, \sigma)$ generated by any two extreme points $\rho,
\sigma$ in $face(\rho, \sigma)$ is a strictly convex, symmetric
and strongly spectral convex set. By Lemma 4.7 $K$ is affinely
isomorphic to the state space of a spin factor, i.e. to the unit
ball in a real Hilbert space. This means that $K$ has "the Hilbert
ball property". Therefore by Theorem \ref{T3.1} $K$ is affinely
and topologically isomorphic to the state space of a
\textit{JB}-algebra with $^\ast$-weak topology. The proof is
complete.
\begin{flushright}
$\square$
\end{flushright}

\section*{Acknowledgments}
\textit{This work was done within
the framework of the Associateship Scheme of the Abdus Salam
International Centre  for Theoretical Physics (ICTP), Trieste,
Italy. The first author would like to thank ICTP for providing financial support
 and all facilities during his visit to ICTP (July-August, 2011).}

\end{document}